\documentclass{article}
\usepackage{amsfonts}

\usepackage{amssymb}

\usepackage{amsmath}

\usepackage{cite}

\newtheorem{theorem}{Theorem}[section]
\newtheorem{definition}[theorem]{Definition}
\newtheorem{lemma}[theorem]{Lemma}
\newtheorem{corollary}[theorem]{Corollary}
\newtheorem{proposition}[theorem]{Proposition}

\date{}

\begin{document}

\title{Free differential Lie Rota-Baxter algebras and Gr\"{o}bner-Shirshov bases\footnote{Supported by the NNSF of China (11171118, 11571121) and  the NSF of Guangdong (2015A030310502). }}

\author{
Jianjun Qiu   \\
{\small \ School of Mathematics and Statistics, Lingnan Normal
University}\\
{\small Zhanjiang 524048, P. R. China}\\
{\small jianjunqiu@126.com}\\
 Yuqun Chen\footnote {Corresponding author.}   \\
{\small \ School of Mathematical Sciences, South China Normal
University}\\
{\small Guangzhou 510631, P. R. China}\\
{\small yqchen@scnu.edu.cn}
}

\maketitle \noindent\textbf{Abstract:} We establish  the  Gr\"{o}bner-Shirshov bases theory for differential  Lie $\Omega$-algebras.  As an application, we  give a linear basis of a    free  differential Lie  Rota-Baxter algebra on a set.

\ \


\ \


%\tableofcontents
\numberwithin{equation}{section}
\section{Introduction}

 Let $k$ be a field  and $\lambda\in k$. A differential
algebra of weight $\lambda$ or  a $\lambda$-differential
algebra (\cite{ri34,  GK08,  ko}) is a   $k$-algebra $(R, \cdot)$ together with a differential  operator (of weight $\lambda$)
 $D:R\rightarrow R $ satisfying
$$
D(x\cdot y)=D(x)\cdot y+x\cdot D(y)+\lambda D(x)\cdot D(y),  x, y\in R.
$$
The differential  algebras  were first  studied by J.F.  Ritt \cite{ri34} and have  developed to be an important branch of   mathematics in  both theory and applications (see for instance  \cite{cohn65, ko, sv03}).

 A Rota-Baxter   algebra of weight $\lambda$ or  $\lambda$-Rota-Baxter
algebra (\cite{Bax60, lguo, Ro69})   is a
 $k$-algebra $(R, \cdot)$ together with a   Rota-Baxter  operator (of weight $\lambda$)
$P:R\rightarrow R$   satisfying
$$
P(x)\cdot P(y)=P(x\cdot P(y))+P(P(x)\cdot y)+ \lambda P(x\cdot y), x, y\in R.
$$
The Rota-Baxter operator on an associative algebra  initially appeared   in  probability    \cite{Bax60} and  then      in combinatorics  \cite{Ro69}  and quantum field theory \cite{ck00}.  There are a number of studies on   associative Rota-Baxter algebras on both commutative and    noncommutative case. For more details we refer the reader to   \cite{lguo} and the references given there.  The  Rota-Baxter operator  of weight 0 on a  Lie algebra  is also   called the operator form of the classical Yang-Baxter equation  \cite{sts}.  The Lie Rota-Baxter algebras  are  closely
related with the  pre-Lie algebras.  Recently,    there are many  results on Lie Rota-Baxter  algebras and   related topics (see for instance
  \cite{ab, bgn10, gp16,  pgb, qc16}).

 Similarly to the relation between  differential operator and integral operator as in the First Fundamental Theorem of Calculus,   L.
Guo and W. Keigher \cite{GK08} introduced the notion of
 differential Rota-Baxter   algebra which is a    $k$-algebra
$R$ together with  a  differential operator $D$ and a Rota-Baxter
operator $P$ such that $DP=Id_{R}$.

 As we known, the free objects   of  various varieties of linear algebras play an important role. Sometimes, it is difficult to give a linear basis of a free algebra, for example, it is an open problem to  find a linear basis of a free Jordan  algebra. A linear basis of the  free differential associative (resp. commutative and  associative) Rota-Baxter algebra on a set   was given   by L. Guo and W. Keigher \cite{GK08}. In this paper,
we    apply the  Gr\"{o}bner-Shirshov bases method to construct  a free  differential Lie  Rota-Baxter algebra. Especially, we   give a linear basis of a free    differential Lie  Rota-Baxter algebra  on a set.

Gr\"{o}bner bases and Gr\"{o}bner-Shirshov bases  have been proved to be very
useful in different branches of mathematics,  which  were invented
independently by A.I. Shirshov \cite{Sh}, H. Hironaka
\cite{Hi64} and B.  Buchberger \cite{Bu70} on  different types of algebras.  For more details on  the Gr\"{o}bner-Shirshov bases and their applications,  see for instance  the surveys   \cite{BC14,  BK03}, the books    \cite{AL, BKu94, CLO,Ei}  and the papers  \cite{bcq,  qc, qc16, DK10, gg15}.

The     $\Omega$-algebra   was
introduced  by A.G. Kurosh \cite{Ku60}.
 A   differential Lie $\Omega$-algebra  over a field $k$   is a differential Lie algebra $L$  with  a set of  multilinear operators  $\Omega$ on $L$. It is easy to see that a differential Lie Rota-Baxter algebra  is a differential Lie $\Omega$-algebra with a single  operator  satisfying  the Rota-Baxter relation.

The paper is organized as follows.
 In Section 2, we review  the Gr\"{o}bner-Shirshov
bases theory  for   differential  associative  $\Omega$-algebras. In Section 3, we firstly   construct a free differential Lie $\Omega$-algebra by the    differential   nonassociative   Lyndon-Shirshov $\Omega$-words, which is a generalization of the classical nonassociative  Lyndon-Shirshov words. Secondly,  we     establish the  Gr\"{o}bner-Shirshov
bases theory   for   differential  Lie   $\Omega$-algebras.   In Section 4,    we obtain a Gr\"{o}bner-Shirshov basis  of a free $\lambda$-differential Lie  Rota-Baxter algebra and then a linear basis of  such an   algebra is obtained  by  the  Composition-Diamond
lemma for differential  Lie   $\Omega$-algebras.

\section{Gr\"{o}bner-Shirshov
bases for  $\lambda$-differential associative   $\Omega$-algebras}

In this section, we briefly  review the  Gr\"{o}bner-Shirshov bases theory
for   $\lambda$-differential associative   $\Omega$-algebras, which can be found in  \cite{qc}.

\subsection{Free     $\lambda$-differential associative   $\Omega$-algebras}

Let $D$ be a $1$-ary operator  and
$$
\Omega:=\bigcup_{m=1}^{\infty}\Omega_{m},
$$
where $\Omega_{m}$ is a set of $m$-ary operators for any $m\geq 1$.  For any set $Y$, we   define the following notations:

$S(Y)$: the set of all nonempty associative words on $Y$.

$Y^*$:    the set of all   associative words on $Y$ including the empty word $1$.

$Y^{**}$: the set of all nonassociative words on $Y$.\\

$\Delta(Y):= \bigcup_{m=0}^{\infty}\{D^m(y)|  y\in Y\},$  where $D^0(y)=y, y\in Y$.\\

$
\Omega(Y):=\bigcup_{m=1}^{\infty}\left\{ \omega^{(m)}(y_1, y_2, \cdots, y_m)|y_i\in Y,  1 \leq i \leq m,  \omega^{(m)}\in \Omega_{m} \right\}.\\
$

Let $X$ be a set.  Define the differential  associative  and   nonassociative  $\Omega$-words on $X$ as follows.  For $n=0$, define
$
 \langle D, \Omega;  X \rangle _{0}=S(\Delta(X)), \  (D, \Omega;  X) _{0}=  (\Delta(X))^{**}.
$
For $n>0$, define
$$
 \langle D, \Omega;  X\rangle_{n} =S(\Delta(X\cup
\Omega(\langle D, \Omega;  X \rangle_{n-1}))),
$$
$$
(D, \Omega;  X) _{n}=  (\Delta(X\cup
\Omega((D, \Omega;  X)_{n-1})))^{**}.
$$
Set
$$
\langle D, \Omega;  X \rangle=\bigcup_{n=0}^{\infty}\langle D, \Omega;  X \rangle_{n}, \ \ (D, \Omega;  X)=\bigcup_{n=0}^{\infty}(D, \Omega;  X) _{n}.
$$
The elements of  $\langle D, \Omega;  X \rangle$ (resp. $(D, \Omega;  X)$) are called   differential  associative  (resp. nonassociative) $\Omega$-words on $X$. A   differential  associative $\Omega$-word  $u$ is called prime if  $u\in \Delta(X\cup \Omega(\langle D, \Omega;  X\rangle))$.

Let $k$ be a field and $\lambda\in k$. A     $\lambda$-differential associative    $\Omega$-algebra  over $k$ is a  $\lambda$-differential associative $k$-algebra  $R$ together  with  a set of  multilinear   operators    $\Omega$ on $R$.

Let
$
DA\langle \Omega; X\rangle=k\langle D, \Omega;  X \rangle
$
be the semigroup
algebra of $\langle D, \Omega;  X \rangle$. Let  $u=u_1u_2\cdots u_t\in \langle D, \Omega;  X \rangle$,  where each $u_i$ is prime. If $t=1$, i.e.   $u=D^i(u')$ for some $i\geq 0, u'\in  X\cup \Omega(\langle D, \Omega;  X\rangle) $, then we define
$
D(u)=D^{i+1}(u').
$
If $t>1$,     then we recursively define
$$
D(u)=D(u_1)(u_2\cdots u_t
)+u_1D(u_2\cdots u_t)+\lambda D(u_1)D(u_2\cdots u_t).
$$
Extend linearly $D$ to $DA\langle \Omega; X\rangle$.  For  any $\omega^{(m)}\in \Omega_m$, define
 $$
 \omega^{(m)}: \langle D, \Omega;  X \rangle^m\rightarrow \langle D, \Omega;  X \rangle, (u_1, u_2, \cdots, u_m)\mapsto  \omega^{(m)}(u_1, u_2, \cdots, u_m)
 $$
and extend it linearly to $DA\langle \Omega; X\rangle^m$.

 \begin{theorem}(\cite{qc})
 $(DA\langle \Omega; X\rangle, D, \Omega)$ is a free   $\lambda$-differential  associative   $\Omega$-algebra on the set $X$.
\end{theorem}

\subsection{Composition-Diamond lemma for \  $\lambda$-differential  \  associative  $\Omega$-algebras }

Let $\star $ is a symbol, which is not in $X$. By a
differential $\star$-$\Omega$-word we mean any expression in  $\langle D, \Omega; X\cup
\{\star\} \rangle$   with only one occurrence of $\star$. The set of all the
differential $\star$-$\Omega$-words on $X$ is denoted by $\langle D, \Omega;  X \rangle^\star$.
Let $\pi$ be a differential $\star$-$\Omega$-word and $s\in DA\langle \Omega; X\rangle$. Then we call
$
\pi|_{s}=\pi|_{\star\mapsto s}
$
a  differential  $s$-word.

%For any $u\in \langle D, \Omega;  X \rangle$, the number
%$
%dep(u):=\mbox{min}\{n|u\in \langle D, \Omega;  X \rangle_n\}
%$
%is called the depth of  $u$.

Let $deg(u)$ be the number
of   all  occurrences of $x\in X$, $\omega\in \Omega$ and $D$ in $u$. If  $u=u_1u_2\cdots u_m$, where $u_i$ is prime, then the breath   of $u$, denoted by  $bre(u)$, is defined to be  the number $m$.
Define
$$
wt(u)=(deg(u),bre(u), u_1, u_2, \cdots, u_m).
$$
Let  $X$  and $\Omega$ be   well-ordered sets and  assume that $\omega >D$ for any  $\omega\in \Omega$. We define the Deg-lex order $>_{_{Dl}}$ on
$\langle D, \Omega;  X \rangle$  as follows.  For any $u=u_1u_2\cdots u_n$  and $ v= v_1v_2\cdots v_m \in \langle D, \Omega;  X \rangle$, where $u_i, v_j$ are prime,  define
$$
u>_{_{Dl}}v \ \mbox{if} \  wt(u)>wt(v)\  \mbox{ lexicographically},
$$
where if $u_i=\omega (u_{i1}, u_{i2}, \cdots, u_{it}), v_i=\theta(v_{i1}, v_{i2}, \cdots, v_{il})$, $\omega, \theta\in \{D\}\cup \Omega $ and $deg(u_i)=deg(v_i)$, then $u_i>_{_{Dl}}v_i$  if
$$
 (\omega , u_{i1}, u_{i2}, \cdots, u_{it} )>(\theta, v_{i1}, v_{i2}, \cdots, v_{il})\  \mbox{lexicographically}.
$$
It is easy to check that $>_{_{Dl}}$ is a well order on
$\langle D, \Omega;  X \rangle$.  For any $0\neq f\in DA\langle \Omega; X\rangle$,  let   $\bar{f}$  be the
 leading term  of $f$ with respect to the order $>_{_{Dl}}$.
Let us denote   $lc(f)$  the coefficient of the leading term $\bar{f}$ of $f$.

For $1\leq t\leq n$, define
$$
I_{n}^{t}=\{(i_1,i_2,\cdots,i_n)\in \{0,
1\}^{n}|i_1+i_2+\cdots+i_n=t\}.
$$

\begin{lemma}\label{p2.2} (\cite{qc}) If  $ u=u_1u_2\cdots u_n\in \langle D, \Omega;  X \rangle$, where each $u_i$ is prime, then
\begin{eqnarray*}
 D(u)&=&\sum_{(i_1,i_2,\cdots,i_n)\in I_{n}^{1}} D^{i_1}(u_1)D^{i_2}(u_2)\cdots D^{i_n}(u_n)\\
 &&+ \sum_{t=2}^n   \sum_{(i_1,i_2,\cdots,i_n)\in I_{n}^{t}}\lambda ^{t-1}D^{i_1}(u_1)D^{i_2}(u_2)\cdots D^{i_n}(u_n).
\end{eqnarray*}
 \end{lemma}

\begin{lemma} \label{l3.3} (\cite{qc}) Let  $ u=u_1u_2\cdots u_n\in \langle D, \Omega;  X \rangle$, where each $u_i$ is prime.

\begin{enumerate}
\item[(a)] If  $\lambda= 0$, then
$
\overline{D^{i}(u)}=D^{i}(u_1) u_2\cdots u_n
$
and   $lc( D^{i}(u) )=1$.
\item[(b)] If  $\lambda\neq 0$, then
$
\overline{D^{i}(u)} =D^{i}(u_1)D^{i}(u_2)\cdots D^{i}(u_n)
$
and $lc( D^{i}(u))=
 \lambda^{(n-1)i}$.

\end{enumerate}
It follows that if $u,v \in \langle D, \Omega;  X \rangle$ and $u>_{_{Dl}}v$, then
$\overline{D(u)}>\overline{D(v)}$.

\end{lemma}

\begin{proposition}\label{l2} (\cite{qc})
For any $u, v \in \langle D, \Omega;  X \rangle,  \pi \in \langle D, \Omega;  X \rangle^{\star}$, if
$u>_{_{Dl}}v$, then $ \overline{\pi|_u}>_{_{Dl}}\overline{\pi|_v} $.
\end{proposition}

If $\overline{\pi|_{s}}=\pi|_{\overline{s}}$, where $s\in DA\langle \Omega; X\rangle$ and $\pi \in \langle D, \Omega;  X \rangle^{\star}$,  then $\pi|_{s}$  is called a  normal
 differential  $s$-word. Note that not each  differential $s$-word is a normal differential  $s$-word,
for example,
if $u= D(x) P(D^2(\star))$
and $s=xy$, where $P\in \Omega, \ x, y\in X$,  then $\pi|_{s}$ is not a   normal differential $s$-word.
However, if we take
$\pi'=D(x_1) P(\star)$,
then $\pi|_{s}=\pi'|_{D^2(s)}$   and $\pi'|_{ D^2(s)}$ is a normal differential $D^2(s)$-word.

\begin{lemma}\label{l3.4} (\cite{qc})
For any differential $s$-word $\pi|_{s}$, there exist $i\geq 0$ and $\pi'$ such that $\pi|_{s}=\pi'|_{D^{i}(s)}$ and
$\pi'|_{D^{i}(s)}$ is  a normal differential $D^{i}(s)$-word.
\end{lemma}

 Let $f, g \in DA\langle \Omega; X\rangle$.   There are two kinds of compositions.

\begin{enumerate}
\item[(i)]If there exists a $w=\overline{D^{i}(f)} a=b \overline{D^{j}(g)}$ for some $a,b\in
\langle D, \Omega;  X \rangle$ such that $bre(w)< bre(\bar{f})+bre(\bar{g})$, then
we call
$$
(f,g)_{w}=lc(D^{i}(f))^{-1}D^{i}(f)a-lc(D^{j}(g))^{-1}bD^{j}(g)
$$
the intersection composition of $f$ and $g$ with respect to the ambiguity  $w$.

\item[(ii)] If there exists a $\pi \in \langle D, \Omega;  X \rangle^{\star}$ such that
$w=\overline{D^{i}(f)}=\pi|_{_{\overline{D^{^{j}}({g})}}}$, where $\pi|_{_{ D^{^{j}}({g})}}$ is a  normal differential $D^{^{j}}({g}) $-word,   then we call
$$
(f,g)_{w}=lc(D^{i}(f))^{-1}D^{i}(f)-lc(D^{j}(g))^{-1}\pi|_{_{D^{^{j}}({g})}}
$$
the inclusion composition of $f$ and $g$ with respect to the ambiguity  $w$.
\end{enumerate}

Let $S$ be a  subset of  $DA\langle \Omega; X\rangle$. Then the
composition $(f,g)_w$ is called trivial modulo $(S,w)$ if
$$
(f,g)_w=\sum\alpha_i\pi_i|_{_{D^{l_i}(s_i)}},
$$
where each $\alpha_i\in k$,  $\pi_i\in \langle D, \Omega;  X \rangle^{\star}$, $s_i\in
S$, $\pi_i|_{_{D^{l_i}(s_i)}}$ is a normal differential  $D^{l_i}(s_i)$-word and
$\pi_i|_{_{\overline{D^{l_i}(s_i)}}}<_{_{Dl}} w$. If this is the case, we write
$$
(f,g)_w\equiv_{ass} 0 \ \ mod (S,w).
$$
In general, for any two polynomials $p$ and $q$, $ p\equiv_{ass} q \ \ mod
(S,w) $ means that $ p-q=\sum\alpha_i\pi_i|_{_{D^{l_i}(s_i)}}, $ where
each $\alpha_i\in k$,  $\pi_i\in \langle D, \Omega;  X \rangle^{\star}$, $s_i\in S$,
$\pi_i|_{_{D^{l_i}(s_i)}}$ is a  normal differential $D^{l_i}(s_i)$-word and
$\pi_i|_{_{\overline{D^{l_i}(s_i)}}}<_{_{Dl}} w$.

A set $S\subset DA\langle \Omega; X\rangle$ is called a Gr\"{o}bner-Shirshov basis  in  $DA\langle \Omega; X\rangle$ if any composition $(f,g)_w$ of $f,g\in S$  is
trivial modulo $(S,w)$.

\begin{theorem}\label{3.6}(\cite{qc}, Composition-Diamond lemma for   differential associative  $\Omega$-algebras) \ \  Let $S$ be
a  subset of  $DA\langle \Omega; X\rangle$, $Id_{DA}(S)$ the ideal of
 $DA\langle \Omega; X\rangle$ generated by $S$ and
 $>_{_{Dl}}$ the Deg-lex order  on $\langle D, \Omega;  X \rangle $ defined as before.  Then the following
statements are equivalent:
 \begin{enumerate}
\item[(i)] $S $ is a Gr\"{o}bner-Shirshov basis in $DA\langle \Omega; X\rangle$.
\item[(ii)] $ f\in Id_{DA}(S)\Rightarrow
\bar{f}=\pi|_{_{\overline{D^{i}(s)}}}$  for some $\pi\in\langle D, \Omega;  X \rangle^\star $, $s\in S$ and $i\geq 0$.
 \item[(iii)]The set
$$
Irr(S)=  \left\{
w\in \langle D, \Omega;  X \rangle
 \left|
 \begin{array}{ll}
  w \neq
\pi|_{_{\overline{D^{i}(s)}}},\ s\in S,\ i\geq 0,\
\pi|_{_{D^{i}(s)}}   \mbox{is  } \\
 \mbox{a normal  differential}\  D^{i}(s)\mbox{-word}
\end{array}
\right. \right\}
$$
is a  linear basis of the   differential associative  $\Omega$-algebra $DA\langle \Omega; X| S\rangle:=DA\langle \Omega; X\rangle/Id_{DA}(S)$.
\end{enumerate}
\end{theorem}

\section{Gr\"{o}bner-Shirshov
bases for $\lambda$-differential Lie  $\Omega$-algebras }

\subsection{Lyndon-Shirshov words}
In this subsection, we review the concept and some properties of Lyndon-Shirshov words, which can be found in \cite{bc07, Sh}.

For any $u\in X^*$, let us denote by  $deg(u)$    the degree (length) of $u$.  Let $>$  be  a   well order on $X$.    Define the lex-order $>_{lex}$ and the deg-lex order $>_{deg-lex}$ on $X^*$ with respect to $>$ by:

 (i) $1>_{lex} u$ for  any nonempty word  $u$, and  if $u=x_iu'$ and $v=x_jv'$, where $x_i, x_j\in X$,   then $u>_{lex} v$ if   $x_i> x_j$, or $x_i=x_j$ and $u'>_{lex} v'$ by induction.

 (ii) $u>_{deg-lex}v$ if $deg(u)>deg(v)$, or $deg(u)=deg(v)$ and $u>_{lex}v$.

A nonempty associative  word $w$ is called  an associative Lyndon-Shirshov word on $X$,  if $  w=uv >_{lex} vu$ for any  decomposition of $w=uv$, where $1\neq  u, v\in X^*$.

A nonassociative word $(u)\in X^{**}$ is said to be   a nonassociative Lyndon-Shirshov word on $X$ with respect to the  lex-order $>_{lex}$,   if
  \begin{itemize}
    \item [(a)] $u$ is an associative Lyndon-Shirshov word on $X$;
    \item [(b)] if $(u)=((v)(w))$, then both $(v)$ and $(w)$ are nonassociative Lyndon-Shirshov words on $X$;
    \item [(c)] if $(v)=((v_1)(v_2))$, then $v_2 \leq_{lex} w$.
  \end{itemize}

Let $ALSW(X)$    (resp. $NLSW(X))$ denote the set of all the associative (resp. nonassociative) Lyndon-Shirshov words on $X$ with respect to the lex-order $>_{lex}$.
It is well known that for any  $u\in ALSW(X)$,   there exists a unique Shirshov standard  bracketing way  $[u]$ (see for instance \cite{bc07}) on $u$ such that $[u] \in NLSW(X)$.
Then
$
NLSW(X )=\{[u] | u\in ALSW(X )\}.
$

Let $k\langle X\rangle$  be the free associative algebra on $X$ over a field $k$ and  $Lie(X)$ be  the Lie subalgebra of   $k\langle X\rangle$  generated by   $X$ under the Lie bracket $(u v)=uv-vu$. It is well known that $Lie(X)$ is a free Lie algebra on the set $X$ and  $NLSW(X)$ is  a linear basis of $Lie(X)$.

\subsection{Differential Lyndon-Shirshov  $\Omega$-words}

 Let  $>_{_{Dl}}$ be the Deg-lex  order on $\langle D, \Omega;  X \rangle$
and $\succ$  the restriction of   $>_{_{Dl}}$ on  $\Delta(X\cup \Omega(\langle D, \Omega;  X \rangle))$. Define the      differential Lyndon-Shirshov  $\Omega$-words  on the set $X$ as follows.

For $n=0$, let $Z_0:=\Delta(X )$. Define
$$
ALSW(D, \Omega;  X)_0:=ALSW(Z_0),
$$
$$
NLSW(D, \Omega;  X)_0:=NLSW(Z_0)= \{[u]|u\in ALSW(D, \Omega;  X)_0\}
$$
with respect to the lex-order $\succ_{lex}$ on $(Z_0)^*$, where $[u]$ is the Shirshov  standard   bracketing way on $u$.

Assume that we have defined
$$
ALSW(D, \Omega;  X)_{n-1},
$$
$$
NLSW(D, \Omega;  X)_{n-1}:= \{[u]|u\in ALSW(D, \Omega;  X)_{n-1}\}.
$$

Let $Z_n:=\Delta(X\cup \Omega(ALSW(D, \Omega;  X)_{n-1}))$.
Define
$$
ALSW(D, \Omega;  X)_n:=ALSW(Z_n).
$$
with respect to the lex-order $\succ_{lex}$ on $Z_n^*$.
For any $u\in Z_n$, define the bracketing way on $u$ by
$$
[u]:= \left\{
 \begin{array}{ll}
u,  & if \  u=D^i(x), x\in X, \\
D^i(\omega^{(m)}([ u_1], [u_2], \cdots,[ u_m])), & if \  u= D^i(\omega^{(m)}(u_1, u_2, \cdots, u_m)).
 \end{array}
  \right.
$$
Let
$
[Z_n]:=\{[u]|u\in Z_n\}.
$
Thus, the order  $\succ$ on $Z_n$ induces an order  on  $[Z_n]$ by $[u]\succ[v]$  if $u\succ v$ for any $u, v\in Z_n$.  For any $u=u_1u_2\cdots u_t\in ALSW(D, \Omega;  X)_n$, where each  $u_i\in Z_n$,   we define
$$
[u]: =[[ u_1 ][ u_2]\cdots [ u_t ]]
$$
the Shirshov standard  bracketing way  on the word $[ u_1 ][ u_2]\cdots [ u_t ]$,  which means that $[u]$  is a nonassociative  Lyndon-Shirshov   word on the set $\{[ u_1 ],[ u_2],\cdots ,[ u_t ]\}$.     Define
$$
NLSW(D, \Omega;  X)_{n}:= \{[u]|u\in ALSW(D, \Omega;  X)_{n}\}.
$$
It is easy to see that
$
NLSW(D, \Omega;  X)_{n}=NLSW([Z_n])
$
with respect to the lex-order $\succ_{lex}$ on $[Z_n]^*$.

Set
$$
ALSW(D, \Omega;  X):=\bigcup_{n=0}^{\infty}ALSW(D, \Omega;  X)_n,
$$
$$
NLSW(D, \Omega;  X):=\bigcup_{n=0}^{\infty}NLSW(D, \Omega;  X)_n.
$$
Then, we have
$$
NLSW(D, \Omega;  X)=\{[u]|u\in  ALSW(D, \Omega;  X)\}.
$$
The elements of $ALSW(D, \Omega;  X)$ (resp. $NLSW(D, \Omega;  X)$) are called the      differential associative (resp. nonassociative) Lyndon-Shirshov $\Omega$-words on the set $X$.

\subsection{Free     $\lambda$-differential Lie  $\Omega$-algebras}

In this subsection,   we prove that the set $NLSW(D, \Omega;  X)$ of all differential  nonassociative  Lyndon-Shirshov   $\Omega$-words on  $X$ forms a linear basis of the  free   $\lambda$-differential Lie   $\Omega$-algebra on  $X$.

A    $\lambda$-differential Lie  algebra is a Lie algebra $L$ with a linear  operator  $D: L\rightarrow L$  satisfying the differential relation
$$
D([xy])=[D(x)y]+[xD(y)]+\lambda [D(x)D(y)],   x, y\in L.
$$
A    $\lambda$-differential Lie   $\Omega$-algebra  is a    $\lambda$-differential  Lie algebra  $L$  with a set of multilinear operators  $\Omega$ on¡¡ $L$.

Let $(R, \cdot, D, \Omega)$ is a    $\lambda$-differential associative  $\Omega$-algebra.  Then it is easy to check that $(R, [, ], D, \Omega)$ is a  $\lambda$-differential Lie  $\Omega$-algebra under   the Lie bracket
$
[a, a']=a\cdot a'-a'\cdot a,\   a, a'\in R.
$

Let $DLie(\Omega; X)$ be the  $\lambda$-differential Lie $\Omega$-subalgebra of  $DA\langle \Omega; X\rangle$  generated by $X$ under the Lie bracket
$
(uv)=uv-vu.
$\\

 Similar to the proofs of Lemma 2.6 and Theorem 2.8 in \cite{qc16}, we have the following results.

\begin{lemma}\label{le2.8}
If  $u\in ALSW(D, \Omega;  X)$, then $\overline{[u]}=u$ with respect to  the order $>_{_{Dl}} $ on $\langle D, \Omega;  X \rangle$.
\end{lemma}

\begin{theorem}
$DLie(\Omega; X)$ is a  free   $\lambda$-differential Lie $\Omega$-algebra on the set $X$ and $NLSW(D, \Omega;  X)$ is a linear basis of $DLie(\Omega; X)$.
\end{theorem}

\subsection{Composition-Diamond lemma for   differential Lie  $\Omega$-algebras }

In this subsection, we  establish the  Composition-Diamond lemma for        differential Lie   $\Omega$-algebras.

\begin{lemma} Let  $\pi\in \langle D, \Omega;  X \rangle^\star$ and $v, \pi|_v\in ALSW(D, \Omega;  X)$. Then there  is a  $\pi'\in \langle D, \Omega;  X \rangle^\star$ and $c\in \langle D, \Omega;  X \rangle$ such that
$$
[\pi|_v]=[\pi'|_{[vc]}],
$$
where $c$ may be empty. Let
$$
[\pi|_v]_{v}=[\pi'|_{[vc]}]|_{[vc]\mapsto [\cdots[[[v][c_1]][c_2]]\cdots [c_m]]}
$$
where   $c=c_1c_2\cdots c_m$ with  each  $c_i\in ALSW(D, \Omega;  X)$ and  $c_t \preceq_{lex} c_{t+1}$.  Then,
$$
[\pi|_v]_{v}=\pi|_{[v]}+\sum \alpha_i \pi_i|_{[v]},
$$
where each $\alpha_i\in k$ and $\pi_i|_{v}<_{_{Dl}} \pi|_{v}$. It follows that $ \overline{[\pi|_v]_{v}}=\pi|_v$   with respect to the order $>_{_{Dl}}$.
\end{lemma}
 {\bf Proof.} The proof is the same as the one of  Lemma 3.2 in  \cite{qc16}. \hfill $\square$\\

Let $0\neq f\in DLie(\Omega; X)\subseteq DA\langle \Omega; X\rangle$.     If $\pi|_{\bar{f}} \in  ALSW(D, \Omega;  X)$, then we call
$$
[\pi|_{f}]_{\bar{f}}=[\pi|_{\bar{f}}]_{\bar{f}}|_{[\bar{f}]\mapsto f}
$$
a special normal differential  $f$-word.

\begin{corollary}\label{co3.10}
Let $f\in DLie(\Omega; X)$ and $\pi|_{\bar{f}}\in ALSW(D, \Omega;  X)$. Then
$$
[\pi|_{f}]_{\bar{f}}=\pi|_{f}+\sum\alpha_i \pi_i|_{f},
$$
\end{corollary}
where each $\alpha_i\in k$ and $\pi_i|_{\bar{f}} <_{_{Dl}}  \pi|_{\bar{f}}$.\\

 Let $f, g \in DLie(\Omega; X)$.  There are two kinds of compositions.

\begin{enumerate}
\item[(i)]If there exists a $w=\overline{D^{i}(f)} a=b \overline{D^{j}(g)}$ for some $a,b\in
\langle D, \Omega;  X \rangle$ such that $bre(w)< bre(\bar{f})+bre(\bar{g})$, then
we call
$$
\langle f,g\rangle_{w}=lc(D^{i}(f))^{-1}[D^{i}(f)a]_{_{\overline{D^{i}(f)}}}-lc(D^{j}(g))^{-1}[bD^{j}(g)]_{_{\overline{D^{j}(g)}}}
$$
the intersection composition of $f$ and $g$ with respect to the ambiguity  $w$.

\item[(ii)] If there exists a $\pi \in \langle D, \Omega;  X \rangle^{\star}$ such that
$w=\overline{D^{i}(f)}=\pi|_{_{\overline{D^{^{j}}({g})}}}$, where $\pi|_{_{ D^{^{j}}({g})}}$ is a  normal differential $D^{^{j}}({g}) $-word,   then we call
$$
\langle f,g\rangle_{w}=lc(D^{i}(f))^{-1}D^{i}(f)-lc(D^{j}(g))^{-1}[\pi|_{_{D^{^{j}}({g})}}]_{_{\overline{D^{^{j}}({g})}}}
$$
the inclusion composition of $f$ and $g$ with respect to the ambiguity  $w$.
\end{enumerate}

If $S$ is a  subset of $DLie(\Omega; X)$, then
the composition $\langle f,g\rangle_w$ is called trivial modulo $(S, w)$ if
$$
\langle f,g\rangle_w=\sum \alpha_i [\pi_i|_{_{D^{l_i}(s_i)}}]_{_{\overline{D^{l_i}(s_i)}}},
$$
where each $\alpha_i\in k, \ s_i\in S$, $[\pi_i|_{_{D^{l_i}(s_i)}}]_{_{\overline{D^{l_i}(s_i)}}}$ is a  special normal differential $D^{l_i}(s_i)$-word   and $\pi_i|_{_{\overline{D^{l_i}(s_i)}}}<_{_{Dl}} w$. If
this is the case, then we write
$$
\langle f,g\rangle_w\equiv 0\  mod(S,w).
$$

In general, for any two polynomials $p$ and $q$, $ p\equiv q \ \ mod
(S,w) $ means that $ p-q=\sum\alpha_i[\pi_i|_{_{D^{l_i}(s_i)}}]_{_{\overline{D^{l_i}(s_i)}}}, $ where
each $\alpha_i\in k$,  $\pi_i\in \langle D, \Omega;  X \rangle^{\star}$, $s_i\in S$,
$[\pi_i|_{_{D^{l_i}(s_i)}}]_{_{\overline{D^{l_i}(s_i)}}}$ is a    normal differential  $D^{l_i}(s_i)$-word and
$\pi_i|_{_{\overline{D^{l_i}(s_i)}}}<_{_{Dl}} w$.

\begin{definition} A   set $S\subset  DLie(\Omega; X)$ is called a
Gr\"{o}bner-Shirshov basis  in   $DLie(\Omega; X)$ if any
composition $\langle f,g\rangle_w$  of    $f, g\in S$ is trivial modulo $(S,w)$.
\end{definition}

\begin{lemma}\label{l3.14}
 Let $f,g\in DLie(\Omega; X)$. Then
$$
\langle f,g\rangle_w-(f,g)_w\equiv_{ass}0 \ \ mod(\{f,g\},w).
$$
\end{lemma}
{\bf Proof.} If $\langle f,g\rangle_w$ and $(f,g)_w$ are compositions of intersection, where $w=\overline{D^i(f)}a=b\overline{D^j(g)}$, then
\begin{eqnarray*}
&&\langle f,g\rangle_w\\
&=&lc(D^{i}(f))^{-1}[D^{i}(f)a]_{_{\overline{D^{i}(f)}}} -lc(D^{j}(g))^{-1} [bD^j(g)]_{_{\overline{D^j(g)}}}\\
&=&lc(D^{i}(f))^{-1}D^{i}(f)b+\sum \alpha_i a_iD^{i}(f)a_i' -lc(D^{j}(g))^{-1}bD^j(g)-\sum \beta_jb_jD^j(g)b_j'\\
&=& (f,g)_w +\sum \alpha_i a_iD^{i}(f)a_i'-  \sum \beta_jb_jD^j(g)b_j' ,
\end{eqnarray*}
where $a_i\overline{D^i(f)}a_i', \  b_j\overline{D^j(g)}b_j'<_{_{Dl}} w$.  It follows that
$$
\langle f,g\rangle_w-(f,g)_w\equiv _{ass} 0 \ mod(\{f,g\}, w).
$$

If $\langle f,g\rangle_w$ and $(f,g)_w$ are compositions of inclusion, where $w=\bar{f}=\pi|_{_{\overline{D^j(g)}}}$, then
$$
\langle f,g\rangle_w=f-lc(D^{j}(g))^{-1}[\pi|_{_{D^j(g)}}]_{_{ \overline{D^j(g)} }}=f-lc(D^{j}(g))^{-1}\pi|_{_{D^j(g)}}-\sum \alpha_i\pi_i|_{_{D^j(g)}},
$$
where $\pi_i|_{_{\overline{D^j(g)}}}<_{_{Dl}}  w$.
It follows that
$$
\langle f,g\rangle_w-(f,g)_w\equiv _{ass} 0 \ mod(\{f,g\}, w).
$$
The proof is complete. \hfill $ \square$\\
\begin{lemma}\label{le3.15}
Let $S\subset DLie(\Omega; X) \subset DA(\Omega; X)$. Then the following two  statements are equivalent:
\begin{enumerate}
\item[(i)]
$S$ is a Gr\"{o}bner-Shirshov basis
in $ DLie(\Omega; X)$,
\item[(ii)] $S$ is a Gr\"{o}bner-Shirshov basis in
$DA\langle \Omega; X\rangle$.
\end{enumerate}
\end{lemma}
{\bf Proof.} $(i)\Longrightarrow(ii)$.
Suppose that $S$ is a Gr\"{o}bner-Shirshov basis in $DLie(\Omega; X)$. Then,
for any composition $\langle f,g\rangle_w$, we have
$$
\langle
f,g\rangle_w=\sum\alpha_i [\pi_i|_{_{D^{l_i}(s_i)}}]_{_{D^{l_i}(s_i)}},
$$
where each $\alpha_i \in k,\  s_i\in S$, $\pi_i|_{_{\overline{D^{l_i}(s_i)}}}<_{_{Dl}} w$ . By Corollary \ref{co3.10}, we have
$$
\langle f,g\rangle_w=\sum\beta_t\pi_t|_{_{D^{l_t}(s_t)}},
$$
where each $\beta_t \in k,\  s_t\in S$, $\pi_t|_{_{\overline{D^{l_t}(s_t)}}}<_{_{Dl}}w$. Therefore,  by Lemma \ref{l3.14}, we  can obtain that
$$
(f,g)_w\equiv_{ass}0 \ mod(S,w).
$$
Thus, $S$ is a Gr\"{o}bner-Shirshov basis in $DA\langle \Omega; X\rangle$.

$(ii)\Longrightarrow(i)$.   Assume that $S$ is a Gr\"{o}bner-Shirshov basis in
$DA\langle  \Omega; X\rangle$. Then, for any composition $\langle
f,g\rangle_w$ in $S$,    we  have  $\langle
f,g\rangle_w\in DLie(\Omega; X)$ and $\langle
f,g\rangle_w\in Id_{DA}(S)$.  By Theorem \ref{3.6},   $\overline{\langle
f,g\rangle_w}=\pi_1|_{_{\overline{D^{i_1}(s_1)}}}\in ALSW(D, \Omega;  X)$.
Let
$$
h_1=\langle
f,g\rangle_w-\alpha_1[\pi_1|_{_{D^{i_1}(s_1)}}]_{_{\overline{D^{i_1}(s_1)}}},
$$
 where $\alpha_1$ is the coefficient of $\overline{\langle
f,g\rangle_w}$.
Then,
$\overline{h_1}<_{_{Dl}} \overline{\langle f,g\rangle_w}$, $h_1\in Id_{DA}(S)$ and $h_1\in DLie(\Omega; X)$. Now, the result follows   from   induction on $\overline{\langle f,g\rangle_w}$.   \hfill $\square$\\

\begin{lemma}\label{le3.16}
Let  $S\subset DLie(\Omega; X)$  and
$$
Irr(S)= \{[w]|w\in ALSW(D, \Omega;  X), w\neq \pi|_{_{\overline{D^i(s)}}},  s\in S,  \pi\in\langle D, \Omega;  X \rangle^\star, i\geq 0\}.
$$
Then, for any  $h\in DLie(\Omega; X)$, $h$ can be expressed by
$$
h=\sum\alpha_i[u_i]+
\sum\beta_j[\pi_j|_{_{D^{l_j}(s_j)}}]_{_{\overline{D^{l_j}(s_j)}}},
$$
where each $\alpha_i, \ \beta_j\in k, \ u_i\in ALSW(D, \Omega;  X), u_i\leq_{Dl}  \bar{h}$ and $s_j\in S$, $\pi_j|_{_{\overline{D^{l_j}(s_j)}}}\leq_{Dl} \bar{h}$.
\end{lemma}
{\bf Proof.} By induction on $\overline{h}$,  we can obtain the result.
 \hfill $ \square$\\

 The following theorem is the   Composition-Diamond lemma  for  differential  Lie $\Omega$-algebras. It is a generalization of  Shirshov's Composition
lemma for Lie algebras \cite{Sh}, which was specialized to
associative algebras by L.A. Bokut \cite{Bo76}, see also G.M.  Bergman
\cite{Be78} and B. Buchberger \cite{bu65, Bu70}.

\begin{theorem}\label{th3.18}
(Composition-Diamond lemma for     differential Lie  $\Omega$-algebras) Let $S\subset  DLie(\Omega; X)$ be a nonempty  set and $Id_{DLie}(S)$ the ideal of $DLie(\Omega; X)$ generated by $S$.
 Then the following statements are equivalent:
\begin{enumerate}
\item[(I)] $S $ is a Gr\"{o}bner-Shirshov basis  in $DLie(\Omega; X)$.
\item[(II)] $f\in Id_{DLie}(S)\Rightarrow \bar{f}=\pi|_{_{  \overline{ D^i(s)}}}\in ALSW(D, \Omega;  X)$ for some $s\in S$, $\pi\in\langle D, \Omega;  X \rangle^\star$ and $i\geq 0$.
\item[(III)] The set
$$
Irr(S)= \{[w]|w\in ALSW(D, \Omega;  X), w\neq \pi|_{_{\overline{ D^i(s)}}}, s\in S,  \pi\in \langle D, \Omega;  X\rangle^\star, i\geq 0 \}
$$
is a  linear basis of  the    $\lambda$-differential Lie  $\Omega$-algebras $DLie(\Omega; X|S)$.
\end{enumerate}
\end{theorem}
{\bf Proof.} $(I)\Longrightarrow(II)$.   Since $f\in Id_{DLie}(S)\subseteq Id_{DA}(S)$, by Lemma  \ref{le3.15} and  Theorem \ref{3.6}, we have
$\bar{f}=\pi|_{\overline{D^i(s)}}$ for some $s\in S$, $\pi\in\langle D, \Omega;  X\rangle^\star$ and $i\geq 0$.

$(II)\Longrightarrow(III)$. \ Suppose that
$
\sum\alpha_i[u_i]=0
$
in $DLie(\Omega; X|S)$, where each $[u_i]\in Irr(S)$ and  $u_i>_{_{Dl}}  u_{i+1}$.  That
is,
$
\sum\alpha_i[u_i]\in{Id_{DLie}(S)}.
$
Then each $\alpha_i$ must be 0. Otherwise, say $\alpha_1\neq0$,
since
$
\overline{\sum\alpha_i[u_i]}=u_1
$
and by (II), we have  $[u_1]\in Irr(S)$,    a contradiction. Therefore, $Irr(S)$ is linear independent. By Lemma \ref{le3.16},  $Irr(S)$
is a  linear basis of     $DLie(\Omega; X|S)=DLie(\Omega; X)/Id_{DLie}(S)$.\\

$(III)\Longrightarrow(I)$. \ For any composition $\langle
f,g\rangle_w$ with $f,g\in{S}$, we have $\langle
f,g\rangle_w\in{Id_{DLie}(S)}$. Then, by (III) and by Lemma \ref{le3.16},
$$
 \langle
f,g\rangle_w=\sum\beta_j[\pi|_{_{D^{l_j}(s_j)}}]_{_{\overline{D^{l_j}(s_j)}}}
$$
where each $\beta_j\in k, \ \pi|_{_{\overline{D^{l_j}(s_j)}}} <_{_{Dl}} w$. This proves that $S$ is a
Gr\"{o}bner-Shirshov basis in $DLie(\Omega; X)$.   \hfill $ \square$\\

\section{Free  $\lambda$-differential Rota-Baxter Lie  algebras}

In this  section, by using Theorem  \ref{th3.18}  we give  a Gr\"{o}bner-Shirshov basis of a     free  $\lambda$-differential Rota-Baxter Lie  algebra on a  set $X$ and then  a linear basis of such an algebra is obtained.

\subsection{Gr\"{o}bner-Shirshov bases for free  $\lambda$-differential Lie  Rota-Baxter algebras}

Let $k$ be a field and $\lambda\in k$. A differential Lie  Rota-Baxter  algebra of weight $\lambda$, called also $\lambda$-differential Lie Rota-Baxter
 algebra,
is a Lie   algebra $L$ with two  linear operators
$P,D:L\rightarrow L$ such that for any $x,y\in L$,
\begin{enumerate}
\item[(a)]
 (Rota-Baxter relation) $[P(x)P(y)]=P([xP(y)])+P([P(x)y])+\lambda
P([xy]);$

\item[(b)](differential relation) $ D([xy])=[D(x)y]+[xD(y)]+\lambda
[D(x)D(y)]$;

\item[(c)] (section relation) $D(P(x))=x$.
\end{enumerate}

It is easy to see that  any  $\lambda$-differential Lie Rota-Baxter
algebra is a     $\lambda$-differential  Lie $\{P\}$-algebra satisfying the relations (a) and (c).

Let $DLie(\{P\}; X)$ be the free  $\lambda$-differential  Lie $\{P\}$-algebra on the set $X$ and write
$$
g(u):=D(P([u]))-[u],
$$
$$
f(u,v):=[P([u])P([v])]-P(([u]P([v])))-P((P([u])[v]))-\lambda P([u][v]),  u>_{_{Dl}}v,
$$
where  $u,  v\in ALSW(D, \{P\};  X)$. Set
$$
S=\{f(u, v), g(w)|u,  v, w\in ALSW(D, \{P\};  X),  u>_{_{Dl}}v\}.
$$

It is clear that
$ DRBL(X):=DLie(\{P\}; X|S)$ is a free $\lambda$-differential Lie Rota-Baxter
algebra on  $X$.

For any $f\in DLie(\{P\}; X)$, let us denote
$r(f):=f-lc(f)[\overline{f}]$.

\begin{lemma} The set
$S_1:=\{D(P([u]))-[u]| u\in ALSW(D, \{P\};  X)\}$ is a   Gr\"{o}bner-Shirshov
basis in $DLie(\{P\}; X)$.
\end{lemma}{\bf Proof.} It is easy to check that  $S_1$   is a Gr\"{o}bner-Shirshov
basis in $DLie(\{P\}; X)$. \hfill $\square$

\begin{lemma} \label{l4.2} Let    $u,  v\in ALSW(D, \{P\};  X)$ and $ u>_{_{Dl}}v$.

\begin{enumerate}
\item[(a)] If $\lambda\neq 0$ and  $j> 0$, then
$$
 D^{j}(f(u,v))
\equiv \lambda^j([\overline{D^{j}(f(u,v))}]- (D^{j-1}([u])D^{j-1}([v])))\  mod(S_1, \overline{D^{j}(f(u, v))}).
$$
\item[(b)]  If $\lambda=0$ and  $j> 0$, then
 $$
 D^{j}(f(u,v))
\equiv [\overline{D^{j}(f(u,v))}]-(D^{j-1}([u])P([v]))\  mod(S_1, \overline{D^{j}(f(u, v))}).
 $$
\end{enumerate}
\end{lemma}
{\bf Proof.}  $(a)$  The proof is by  induction on $j$.
  For $j=1$,  we have
\begin{eqnarray*}
&& D(f(u,v)) \\
&=& D((P([u])P([v])))-D(P(([u]P([v])))-D(P((P([u])[v])))-\lambda D(P(([u][v])))\\
&\equiv&\lambda
(D(P([u]))D(P([v])))-\lambda D(P(([u][v])))\\
&\equiv& \lambda ([\overline{D(f(u,v))}]-([u][v]))\ mod(S_1,\overline{D(f(u,v))}).
\end{eqnarray*}

Assume that the result is true for $j-1, j\geq 2$, i.e.
\begin{eqnarray*}
D^{j-1}(f(u,v))
&=&\lambda^{j-1}
(D^{j-1}(P([u]))D^{j-1}(P([v])))\\
&&-\lambda^{j-1}(D^{j-2}([u])D^{j-2}([v]))
+\sum \alpha_i
[\pi_i|_{_{D^{t_i}(s_i)}}]_{_{\overline{D^{t_i}(s_i)}}},
\end{eqnarray*}
where each $\alpha_i\in k$, $s_i\in
S_1$, $\pi_i|_{_{\overline{D^{t_i}(s_i)}}}<_{_{Dl}}D^{j-1}(P(u))D^{j-1}(P(v))$. Since $S_1$ is a Gr\"{o}bner-Shirshov
basis in $DLie(\{P\}; X)$,
$$
D(\sum \alpha_i [\pi_i|_{_{D^{t_i}(s_i)}}]_{_{\overline{D^{t_i}(s_i)}}})=\sum
\beta_{l}[\sigma_l|_{_{D^{n_l}(s_l)}}]_{_{\overline{D^{n_l}(s_l)}}},
$$
where each $\beta_l\in k, s_l\in S_1$, $[\sigma_l|_{_{D^{n_l}(s_l)}}]_{_{\overline{D^{n_l}(s_l)}}}$ is a special  normal  differential $D^{k_l}(s_l)$-word. By Lemma \ref{l3.3},
\begin{eqnarray*}
\overline{[\sigma_l|_{_{D^{n_l}(s_l)}}]_{_{\overline{D^{n_l}(s_l)}}}}&=&\sigma_l|_{_{ \overline{D^{n_l}(s_l)}}}\\
&<_{_{Dl}}&
\overline{D((D^{j-1}((P([u]))D^{j-1}(P([v]))))}=D^{j}(P(u))D^{j}(P(v)).
\end{eqnarray*}
Thus,    we have
\begin{eqnarray*}
&&D^{j}(f(u,v))\\
&=&D(D^{j-1}(f(u,v)))\\
&\equiv& \lambda^{j-1}D((
D^{j-1}(P([u]))D^{j-1}(P([v]))))-\lambda^{j-1}D((D^{j-2}([u])D^{j-2}([v])))\\
&\equiv& \lambda^{j}(D^{j}(P([u]))D^{j}(P([v])))+  \lambda^{j-1}(D^{j}(P([u]))D^{j-1}(P([v])))\\
&&+ \lambda^{j-1}(D^{j-1}(P([u]))D^{j}(P([v])))- \lambda^{j} (D^{j-1}([u])D^{j-1}([v]))\\
&&-  \lambda^{j-1} (D^{j-1}([u])D^{j-2}([v]))-\lambda^{j-1} (D^{j-2}([u])D^{j-1}([v]) )\\
&\equiv& \lambda^{j} (D^{j}(P([u]))D^{j}(P([v])))+ \lambda^{j-1} ( D^{j-1}([u])D^{j-2}([v]))\\
&&+ \lambda^{j-1} (D^{j-2}([u])D^{j-1}([v])))
 - \lambda^{j} (D^{j-1}([u])D^{j-1}([v]))\\
 &&-\lambda^{j-1} (D^{j-1}([u])D^{j-2}([v]))-\lambda^{j-1} (D^{j-2}([u])D^{j-1}([v])))\\
&\equiv& \lambda^{j}(D^{j}(P([u]))D^{j}(P([v]))- \lambda^{j} (D^{j-1}([u])D^{j-1}([v])) \\
&\equiv& \lambda^{j}([\overline{D^{j}(f(u, v))}]-  (D^{j-1}([u])D^{j-1}([v])))\   mod(S_1, \overline{D^{j}(f(u, v))}).
\end{eqnarray*}

$(b)$  The proof is similar to   Case (a).
  \hfill $\square$

\begin{theorem}\label{t4.1} With  the order $>_{_{Dl}}$
 on $\langle D, \{P\};  X \rangle$ defined as before,  the set $S$ is a Gr\"{o}bner-Shirshov
basis in $DLie(\{P\}; X)$.
\end{theorem}
\noindent {\bf Proof.} There are two cases $\lambda
\neq 0$ and $\lambda =0$ to consider.

Case 1.  For $\lambda \neq 0$,  all possible
compositions of the polynomials  in $S$
 are list as below:
 \begin{tabbing}
 $\langle g(\pi|_{_{D^{^{i}}(D(P(v)))}}), g(v)\rangle_{w_1},\ \   w_1=D^{j}(D(P(\pi|_{_{D^{^{i}}(D(P(v)))}}))),$\\[0.7ex]
$\langle  g(\pi|_{_{D^{^{i}}(P(u))D^{^{i}}(P(v))}}),  f(u, v)\rangle_{w_2}, \ \ \  w_2=D^{j}(D(P(\pi|_{_{D^{^{i}}(P(u))D^{^{i}}(P(v))}})))$, \\[0.7ex]
$\langle  f(u, v),  g(v)\rangle_{w_3}, \ \   w_3=D^{l}(P(u))D^l(P(v)),  \ l>0,$\\[0.7ex]
$\langle  f(u, v),  g(u)\rangle_{w_4}, \ \    w_4=D^{l}(P(u))D^l(P(v)), \ l>0,$\\[0.7ex]
$\langle  f(\pi|_{_{D^{^{i}}(D(P(u)))}}, v),  g(u)\rangle_{w_5}, \ \    w_5=D^{j}(P(\pi|_{_{D^{^{i}}(D(P(u)))}}))D^{j}(P(v)), $\\[0.7ex]
$\langle  f(u, \pi|_{_{D^{^{i}}(D(P(v)))}})),  g(v)\rangle_{w_6}, \ \    w_6=D^{j}(P(u))D^{j}(P(\pi|_{_{D^{^{i}}(D(P(v)))}})), $\\[0.7ex]
$\langle  f(u, v),  f(v, w)\rangle_{w_7}, \ \    w_7=D^{^{j}}(P(u))D^{^{j}}(P(v))D^{^{j}}(P(w)), $\\[0.7ex]
$\langle  f(\pi|_{_{D^{i}(P(u))D^{i}(P(v))}} , w),  f(u, v)\rangle_{w_8}, \ \    w_8=D^{j}(P(\pi|_{_{D^{i}(P(u))D^{i}(P(w))}}))D^{j}(P(w)),  $\\[0.7ex]
$\langle  f( u, \pi|_{_{D^{i}(P(v))D^{i}(P(w))}}),  f(v, w)\rangle_{w_9},   \    w_9=D^{j}(P(u))D^{j}(P(\pi|_{_{D^{i}(P(v))D^{i}(P(w))}})),$\\[0.7ex]
where $i, j \geq 0$.
\end{tabbing}

We   check   that all the compositions in $S$ are trivial. Here, we just check one composition  as example.

If $j>0$, then by Lemma \ref{l4.2}, we have
\begin{eqnarray*}
&&\langle  f(\pi|_{_{D^{i}(P(u))D^{i}(P(v))}}, w),  f(u, v)\rangle_{w_8}  \\
&=& \lambda^{-j}D^{j}(f(\pi|_{_{D^{i}(P(u))D^{i}(P(v))}}, w))- \lambda^{-i}[D^{j}(P(\pi|_{_{D^{i}(f(u,v))}})D^{j}(P(w))]_{_{\overline{D^(i)(f(u,v))}}}\\
&\equiv&  -(D^{j-1}([\pi|_{_{D^{i}(P(u))D^{i}(P(v))}}])D^{j-1}([w]))-\lambda^{-i}(D^{j}(P(r([\pi|_{_{D^{i}(f(u,v))}}]_{_{\overline{D^{i}(f(u,v))}}})))D^{j}(P([w])))\\
&\equiv& \lambda^{-i} (D^{j-1}(r([\pi|_{_{D^{i}(f(u,v))}}]_{_{\overline{D^{i}(f(u,v))}}}))D^{j-1}([w]))-\lambda^{-i} (D^{j-1}(r([\pi|_{_{D^{i}(f(u,v))}}]_{_{\overline{D^{i}(f(u,v))}}}))D^{j-1}([w]))\\
 &\equiv&0 \ mod(S, w_8).
\end{eqnarray*}

If $j=0$, then
\begin{eqnarray*}
&&\langle  f(\pi|_{_{D^{i}(P(u))D^{i}(P(v))}}, w),  f(u, v)\rangle_{w_8} \\
&=&  f(\pi|_{_{D^{i}(P(u))D^{i}(P(v))}}, w)- \lambda^{-i}[ P(\pi|_{_{D^{i}(f(u,v))}}) P(w)]_{_{\overline{D^(i)(f(u,v))}}}\\
&\equiv&  -P((P([\pi|_{_{D^{i}(P([u]))D^{i}(P([v]))}}])[w]))-P(([\pi|_{_{D^{i}(P(u))D^{i}(P(v))}}]P([w] )))\\
&&-\lambda P(([\pi|_{D^{i}(P(u))D^{i}(P(v))}][w]))
 -\lambda^{-i}(P(r([\pi|_{_{D^{i}(f(u,v))}}]_{_{\overline{D^{i}(f(u,v))}}})) P([w]) )\\
 &\equiv&   \lambda^{-i}P((P(r([\pi|_{_{D^{i}(f(u,v))}}]_{_{\overline{D^{i}(f(u,v))}}}))[w]))+ \lambda^{-i}P(( r([\pi|_{_{D^{i}(f(u,v))}}]_{_{\overline{D^{i}(f(u,v))}}}) )))\\
&& \lambda^{-i+1}  P((r([\pi|_{_{D^{i}(f(u,v))}}]_{_{\overline{D^{i}(f(u,v))}}})[w]))
 -\lambda^{-i}P((P(r([\pi|_{D^{i}(f(u,v))}]_{\overline{D^{i}(f(u,v))}}))[w]))\\
 &&- \lambda^{-i}P(( r([\pi|_{_{D^{i}(f(u,v))}}]_{_{\overline{D^{i}(f(u,v))}}}) )))-\lambda^{-i+1}  P((r([\pi|_{_{D^{i}(f(u,v))}}]_{_{\overline{D^{i}(f(u,v))}}})[w]))\\
&\equiv&0 \ mod(S, w_8).
\end{eqnarray*}

Case 2. For $\lambda = 0$,  all possible
compositions of the polynomials  in $S$ are list as below:
 \begin{tabbing}
$\langle g(\pi|_{_{D^{^{i}}(D(P(v)))}}), g(v)\rangle_{w_1},  \ {w_1}=D^{j}(D(P(\pi|_{_{D^{^{i}}(D(P(v)))}}))),  $\\[0.7ex]
$\langle g(\pi|_{_{D^{^{i}}(P(u))P(v)}}), f(u, v)\rangle_{w_2}, \  {w_2}= D^{j}(D(P(\pi|_{_{D^{^{i}}(P(u))P(v)}}))),$\\[0.7ex]
$\langle f(u, v), g(u)\rangle_{w_3}, \  {w_3}=D^l(P(u))P(v), \  l>0,$\\[0.7ex]
$\langle f(\pi|_{_{D^{^{i}}(D(P(u)))}}, v), g(u)\rangle_{w_4}, \ {w_4}= D^j(P(\pi|_{_{D^{^{i}}(D(P(u)))}}))P(v),$ \\[0.7ex]
$\langle f(u, \pi|_{_{D^{^{i}}(D(P(v)))}}), g(v)\rangle_{w_5}, \ {w_5}= D^j(P(u))P(\pi|_{_{D^{^{i}}(D(P(v)))}}),$ \\[0.7ex]
$ \langle f(u, v), f(v, w) \rangle_{w_6}, \ {w_6}= D^{^{j}}(P(u))P(v)P(w),$  \\[0.7ex]
$ \langle f(\pi|_{_{D^{^{i}}( P(u) )P(v)}} , w), f(u, v) \rangle_{w_7}, \ {w_7}= D^{^{j}}(P(\pi|_{_{D^{^{i}}( P(u) )P(v)}}))P(w),$\\[0.7ex]
$ \langle f(u, \pi|_{_{D^{^{i}}( P(v) )P(w)}}), f(v, w) \rangle_{w_8}, \ {w_8}= D^{^{j}}(P(u))P(\pi|_{_{D^{^{i}}( P(v) )P(w)}})),$
\end{tabbing}
where $i,j\geq 0$.
We check  that all the compositions in $S$ are trivial. The proof is similar to   Case 1.
\hfill $\square$

\subsection{A linear basis  of a free  $\lambda$-differential Lie Rota-Baxter algebra}
In this subsection, by  Theorems \ref{th3.18} and    \ref{t4.1},  we obtain   a linear
basis of the  free $\lambda$-differential Lie Rota-Baxter algebra on the set $X$.

For $n=0$, define $\langle \{P\}; \Delta(X)\rangle _{0}:=S(\Delta(X))$
 and $  (\{P\}; \Delta(X)) _{0}:=  (\Delta(X))^{**}.$
For $n>0$, define
$$
 \langle  \{P\}; \Delta(X) \rangle_{n} :=S(\Delta(X)\cup
P( \langle  \{P\}; \Delta(X) \rangle_{n-1})),
$$
$$
(\{P\}; \Delta(X)) _{n}:=  (\Delta(X)\cup
P(\{P\}; \Delta(X)) _{n-1})^{**}.
$$
Set
$$
\langle  \{P\}; \Delta(X)\rangle:=\bigcup_{n=0}^{\infty}\langle  \{P\}; \Delta(X) \rangle_{n}, \ \ (\{P\}; \Delta(X)):=\bigcup_{n=0}^{\infty}(\{P\}; \Delta(X))_{n}.
$$

Let $\star $ is a symbol, which is not in $X$. By a
$\star$-$\Omega$-word on $\Delta(X )$, we mean any expression in $\langle \{P\}; \Delta(X )\cup
\{\star\} \rangle$ with only one occurrence of $\star$. We will denote by   $\langle \{P\}; \Delta(X )\rangle^\star$   the set of all the
  $\star$-$\Omega$-words on $\Delta(X )$.
Let $\pi$ be a $\star$-$\Omega$-word on $\Delta(X )$ and $u\in \langle \{P\}; \Delta(X)\rangle$.  Let us denote
$
\pi|_{u}=\pi|_{\star\mapsto u},
$
i.e.   $\star$  is replaced by $u$.

 It is easy to see that
$\langle \{P\}; \Delta(X)\rangle\subseteq \langle D, \{P\};  X\rangle$.  We also use the order $>_{_{Dl}}$ on $\langle \{P\}; \Delta(X)\rangle$ and $\succ$ on $\Delta(X)\cup P(\langle \{P\}; \Delta(X)\rangle)$.

For $n=0$, let $Y_0=\Delta(X)$. Define
$$
ALSW(\{P\};  \Delta(X))_0:=ALSW(Y_0),
$$
$$
NLSW(\{P\};  \Delta(X))_0:=NLSW(Y_0)= \{[u]|u\in ALSW(\{P\};  \Delta(X))_0\}
$$
with respect to the lex-order $\succ_{lex}$.

Assume that we have defined
$$
ALSW(\{P\};  \Delta(X))_{n-1},
$$
$$
NLSW(\{P\};  \Delta(X))_{n-1}= \{[u]|u\in ALSW(\{P\};  \Delta(X))_{n-1}\}.
$$

Let $Y_n:=\Delta(X)  \cup P(ALSW(\{P\};  \Delta(X))_{n-1})$.
Define
$$
ALSW(\{P\};  \Delta(X))_n:=ALSW(Y_n).
$$
For  any    $u\in Y_n$, define the bracketing way on $u$ by
$$
[u]:= \left\{
 \begin{array}{ll}
u,  & if \  u\in \Delta(X), \\
P([u_1]),  & if \  u= P(u_1).
 \end{array}
  \right.
$$
Let
$
[Y_n]:=\{[u]|u\in Y_n\}.
$
Therefore,  the order  $\succ$ on $Y_n$ induces an order   on  $[Y_n]$    by $[ u] \succ[ v ]$  if $u\succ v$ for any $u, v\in Y_n$.
For any  $u=u_1u_2\cdots u_m\in ALSW(\{P\};  \Delta(X))_n$, where each  $u_i\in Y_n$, let us denote
$$
[u]: =[[ u_1 ][ u_2]\cdots [ u_m ]]
$$
the   nonassociative  Lyndon-Shirshov   word on  $\{[ u_1 ],[ u_2],\cdots ,[ u_m ]\}$   with respect to the lex-order $\succ_{lex}$.

Define
$$
NLSW(\{P\};  \Delta(X))_n:=\{[  u ]|u\in ALSW(\{P\};  \Delta(X))_n\}.
$$
It is easy to see that $NLSW(\{P\};  \Delta(X))_n=NLSW([Y_n])$.  Define
$$
ALSW(\{P\};  \Delta(X)):=\bigcup_{n=0}^{\infty}ALSW(\{P\};  \Delta(X))_n,$$
$$
NLSW(\{P\};  \Delta(X)):=\bigcup_{n=0}^{\infty} NLSW(\{P\};  \Delta(X))_n.
$$
Therefore,
$$
  NLSW(\{P\};  \Delta(X)) =\{[  u ]|u\in  ALSW(\{P\};  \Delta(X))\}.
$$

By   Theorems \ref{th3.18} and \ref{t4.1},    we have the following
theorem.

\begin{theorem} The set
$$
Irr(S)=  \left\{
[w]\in NLSW (\{P\};  \Delta(X))
 \left|
 \begin{array}{ll}
  w\neq \pi|_{_{P(u)P(v)}},  \ \    \pi\in   \langle \{P\};  \Delta(X) \rangle^\star \\
 u, v \in ALSW(\{P\};  \Delta(X)), u >_{_{Dl}} v  \\
\end{array}
\right. \right\}
$$
is a linear basis of the free    $\lambda$-differential Lie  Rota-Baxter  algebra $DRBL(X)$ on  $X$.
\end{theorem}

%\noindent{\bf Acknowledgement}: The authors would like to thank
%Professor L.A. Bokut for his guidance, useful discussions and
%enthusiastic encouragement in writing up this paper.

\end{document}